\documentclass[12pt]{article}

\usepackage{amsfonts,graphicx}

\newtheorem{lemma}{Lemma}

\newtheorem{corollary}{Corollary}

\title{
Plane geometry and convexity\\ of polynomial stability regions$^1$
}

\author{Didier Henrion$^{2,3}$, Michael \v Sebek$^3$}

\graphicspath{{/home/henrion/Matlab/Geometry/}}

\begin{document}

\maketitle

\footnotetext[1]{This work was partly supported
by project No. MSM6840770038 of the Ministry of Education of the Czech Republic.}
\footnotetext[2]{LAAS-CNRS, University of Toulouse, France}
\footnotetext[3]{Faculty of Electrical Engineering,
Czech Technical University in Prague, Czech Republic}

\begin{abstract}
The set of controllers stabilizing a linear system
is generally non-convex
in the parameter space. In the case of
two-parameter controller design (e.g. PI control
or static output feedback with one input and two outputs),
we observe however that quite often for benchmark problem
instances, the set of stabilizing controllers seems to be
convex. In this note we use elementary techniques from real
algebraic geometry (resultants and B\'ezoutian matrices)
to explain this phenomenon. As a byproduct, we 
derive a convex linear matrix inequality (LMI) formulation
of two-parameter fixed-order controller design problem,
when possible.
\end{abstract}

\begin{center}
{\bf\small Keywords}\\[1em]
control theory; convexity; resultants
\end{center}

\section{Introduction}

Despite its elementary formulation, the problem of fixed-order
controller design for linear time-invariant systems remains
mostly open. Especially scarce are numerically efficient
computer-aided control system design algorithms in the
fixed-order case, sharply contrasting with the large
number of tools available to solve static state feedback
design or dynamical output feedback design with controllers
of the same order as the plant. Mathematically, fixed-order
controller design can be formulated as a non-convex
non-smooth optimization problem in the parameter space. To the
best of our knowledge, randomized algorithms are amongst
the most efficient numerical methods to cope with this
class of difficult problems. See \cite{hifoo,rct} for
computer experiments supporting this claim, thanks to
public-domain Matlab packages (HIFOO and the Randomized
Control System Toolbox).

This note was motivated by the observation, made by the
first author during a workshop at AIM in August 2005 
\cite{aim}, that 6 out of the 7 two-dimensional instances
of static output feedback (SOF) design problems found
in the database COMPleib \cite{compleib}
seem to be convex. Further motivation was provided
by the excellent historical survey \cite{gryazina}
on D-decomposition techniques, previously studied in
deep detail in \cite{siljak} and \cite{ackermann}.
In \cite{gryazina} the authors describe the intricate geometry
of two-dimensional stability regions with the help of
illustrative examples.
Quite often, the stability regions represented
in these references seem to be convex.

In this note, we use basic results from real algebraic
geometry to detect convexity of the stability region
in the two-parameter case (including PI controllers,
PID with constant gain, SOF design
with one input two outputs or two inputs one output).
We also derive, when possible, a linear matrix
inequality (LMI) formulation of the stability region.

\section{Problem statement}

We consider a parametrized polynomial
\begin{equation}\label{poly}
p(s,k) = p_0(s) + k_1 p_1(s) + k_2 p_2(s)
\end{equation}
where the $p_i(s) \in {\mathbb R}[s]$ are given
polynomials of $s \in {\mathbb C}$ and the $k_i \in
{\mathbb R}$ are parameters. We assume, without loss
of generality, that the ratio $p_1(s)/p_2(s)$ is not
a constant.

Define the stability region
\[
{\mathcal S} = \{k \in {\mathbb R}^2 \: :\:
p(s,k) \:\:\mathrm{stable}\}
\]
where stability is meant in the continuous-time
sense, i.e. all the roots of $p(s,k)$ must lie in the
open left half-plane\footnote{Similar results can be derived
for discrete-time unit disk stability or other semialgebraic
stability domains of the complex plane, but this is
not covered here.}.

We are interested in the following problems:
\begin{itemize}
\item Is stability region $\mathcal S$ convex ?
\item If it is convex, give an LMI representation
\[
{\mathcal S} = \{k \in {\mathbb R}^2 \: :\:
F_0 + F_1 k_1 + F_2 k_2 \succ 0\}
\]
when possible, where the $F_i$ are real symmetric
matrices to be found, and $\succ 0$ means positive definite.
\end{itemize}

\subsection{Example: PI controller design}

Let $b(s)/a(s)$ with $a(s),b(s) \in {\mathbb R}[s]$ denote
the transfer function of an open loop plant, and consider
a proportional integral (PI) controller $k_1/s+k_2$ in a 
standard negative feedback configuration. The closed-loop
characteristic polynomial (\ref{poly}) is then $p(s,k) = sa(s)
+(k_1+k_2s)b(s)$ hence $p_0(s) = sa(s)$, $p_1(s)=b(s)$
and $p_2(s)=sb(s)$.

\subsection{Example: static output feedback}

Given matrices $A \in {\mathbb R}^{n \times n}$, $B \in {\mathbb
R}^{n \times m}$, $C \in {\mathbb R}^{p \times n}$, we want to
find a matrix $K \in {\mathbb R}^{m \times p}$ such that the
closed-loop matrix $A+BKC$ is stable. When $mp=2$, the
characteristic polynomial (\ref{poly})
writes $p(s,k) = \det (sI_n-A-BKC)$
hence $p_0(s) = \det (sI_n-A)$, $p_1(s) = \det (sI_n-B[k_1,0]C)$
and $p_2(s) = \det (sI_n-B[0,k_2]C)$.

\section{Hermite matrix}

The Routh-Hurwitz criterion for stability of polynomials
has a symmetric version called the Hermite criterion.
A polynomial is stable if and only if
its Hermite matrix, quadratic in the polynomial coefficients,
is positive definite. In control systems terminology,
the Hermite matrix is a particular choice of a Lyapunov matrix certifying
stability \cite{parks}.
Algebraically, the Hermite matrix can be defined
via the B\'ezoutian, a symmetric form of the resultant
\cite[Section 5.1.2]{mourrain}.

Let $a(u)$, $b(u)$ be two polynomials of degree $n$
of the indeterminate $u$.
Define the B\'ezoutian matrix $B_u(a,b)$ as the symmetric
matrix of size $n$ with entries $b_{ij}$ satisfying
the linear equations
\[
\frac{a(u)b(v)-a(v)b(u)}{v-u} = \sum_{i=1}^n \sum_{j=1}^n
b_{ij} u^{i-1} v^{j-1}.
\]
The polynomial $r_u(a,b) = \det B_u(a,b)$ is the resultant
of $a(u)$ and $b(u)$ with respect to $u$. It is obtained
by eliminating $u$ from the system of equations
$a(u)=b(u)=0$.

The Hermite matrix of $p(s)$ is defined as the B\'ezoutian
matrix of the real part and the imaginary part of
$p(j\omega)$:
\[
\begin{array}{rcl}
p_R(\omega^2) & = & \mathrm{Re}\:p(j\omega) \\
\omega p_I(\omega^2) & = & \mathrm{Im}\:p(j\omega)
\end{array}
\]
that is, $H(p) = B_{\omega}(p_R(\omega^2),\omega p_I(\omega^2))$.
Let us assume that $p(s)$ is monic, with unit leading coefficient.
The Hermite stability criterion can be formulated as follows.

\begin{lemma}\label{hermite}
Polynomial $p(s)$ is stable if and only if $H(p) \succ 0$.
\end{lemma}

{\bf Proof}:
The proof of this result can be found in \cite{parks}
for example. It can also be proved via
Cauchy indices and Hermite quadratic forms for counting
real roots of polynomials, see \cite[Section 9.3]{bpr}
$\Box$.

By construction,
the Hermite matrix of parametrized polynomial (\ref{poly})
\[
H(p(s,k)) = H(k) = \sum_{i_1,i_2=0}^2 H_{i_1 i_2} k_1^{i_1}k_2^{i_2}
\]
is quadratic in $k$. Therefore,
the Hermite criterion yields a quadratic matrix inequality
formulation of the stability region:
\[
{\mathcal S} = \{k \in {\mathbb R}^2 \: :\: H(k) \succ 0\}.
\]
Quadratic matrix inequalities, a generalization
of bilinear matrix inequalities, typically generate
non-convex regions. For example, the scalar quadratic
inequality $k_1^2-1>0$ models a disconnected,
hence non-convex set.

Surprisingly, it turns out that $\mathcal S$,
even though modeled by a quadratic matrix
inequality, is often a convex set
for practical problem instances.
Here are some examples.

\subsection{Examples: static output feedback}

\begin{figure}[h]
\hfill
\begin{minipage}[t]{.3\textwidth}
\begin{center}  
\includegraphics[scale=0.3]{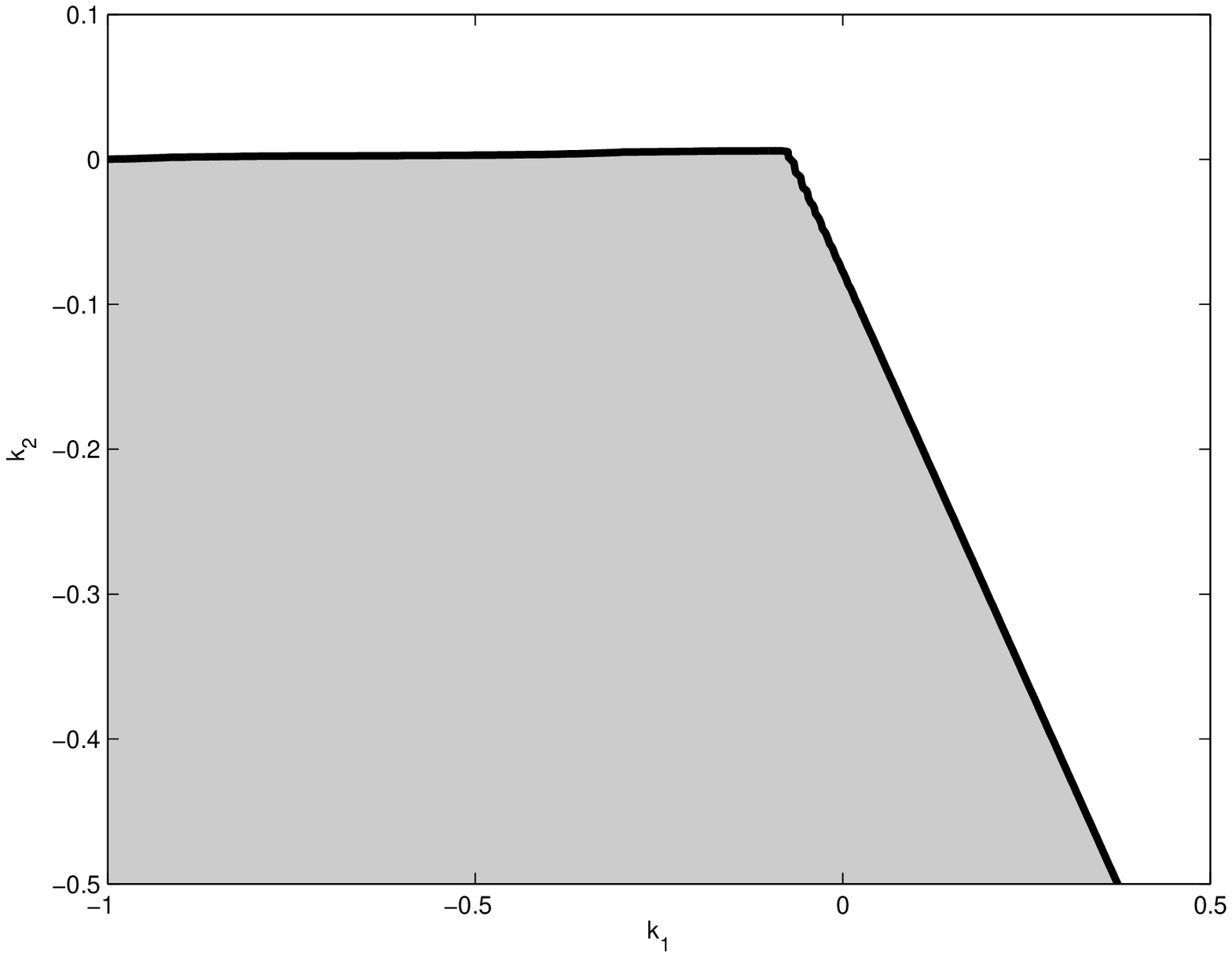}
\caption{{\tt AC4}.\label{fig-ac4}}
\end{center}
\end{minipage}
\hfill
\begin{minipage}[t]{.3\textwidth}
\begin{center}
\includegraphics[scale=0.3]{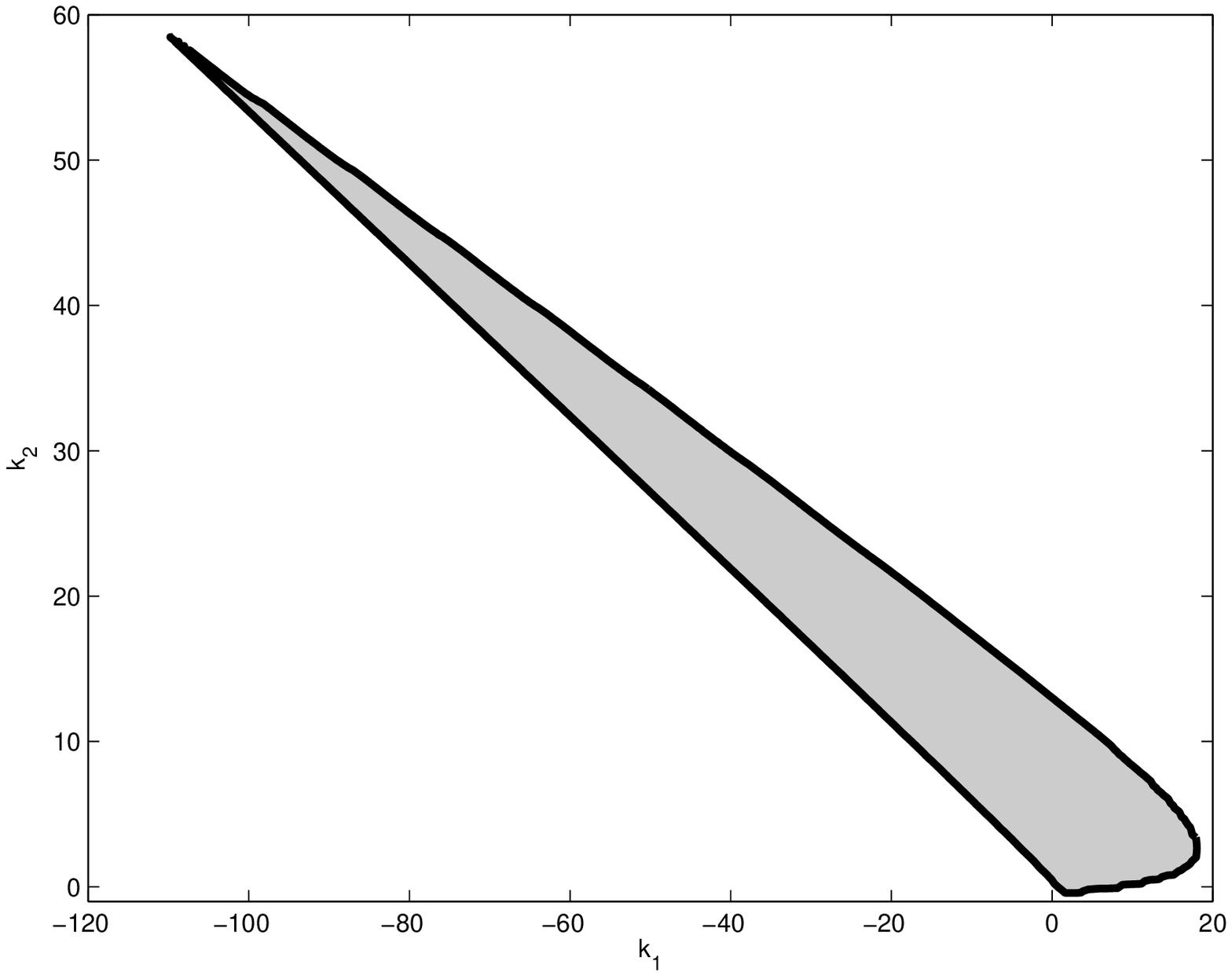}
\caption{{\tt AC7}.\label{fig-ac7}} 
\end{center}
\end{minipage}
\hfill
\begin{minipage}[t]{.3\textwidth}
\begin{center}
\includegraphics[scale=0.3]{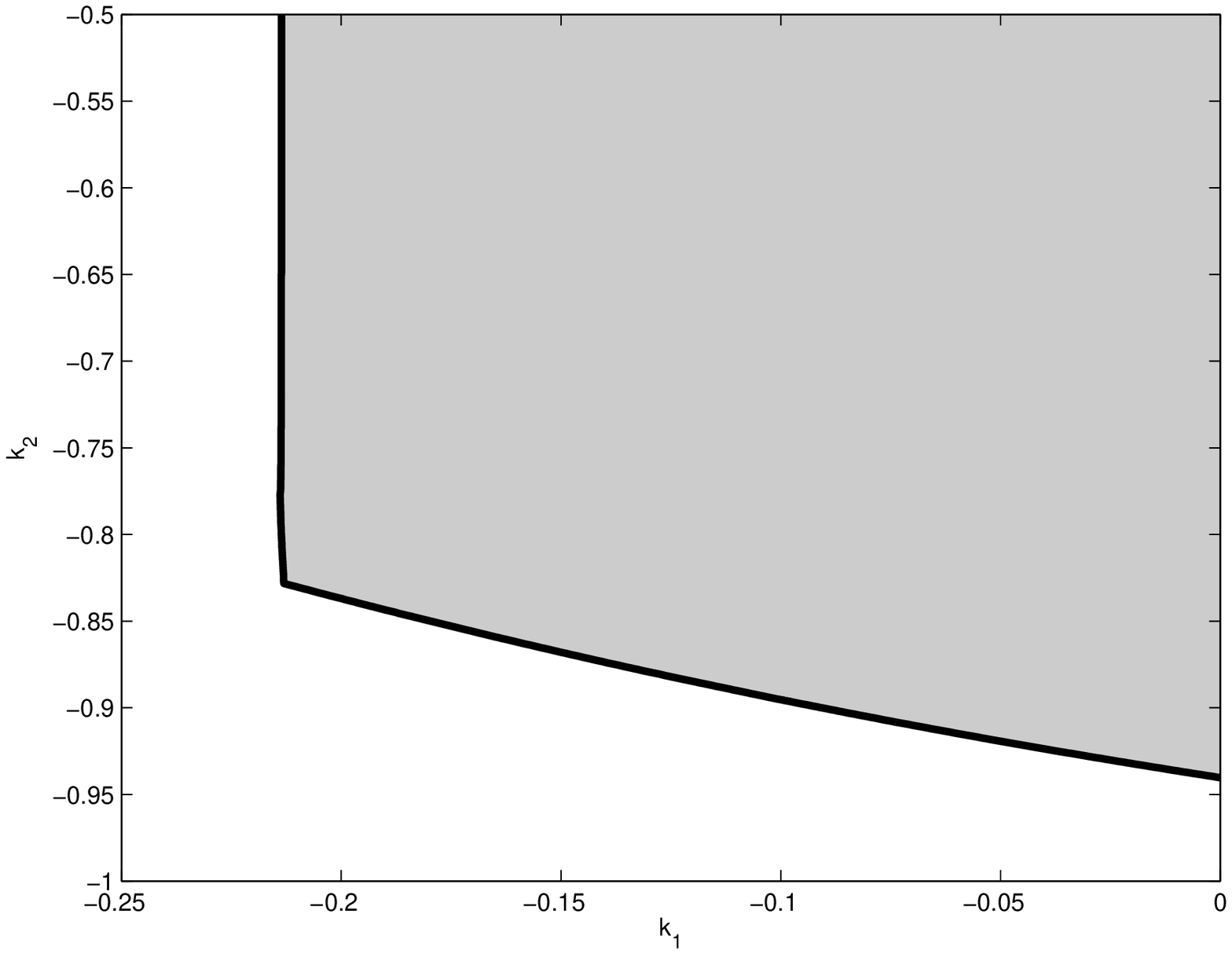}
\caption{{\tt AC17}.\label{fig-ac17}}
\end{center}
\end{minipage}
\end{figure}

\begin{figure}[h]
\hfill
\begin{minipage}[t]{.3\textwidth}
\begin{center}  
\includegraphics[scale=0.3]{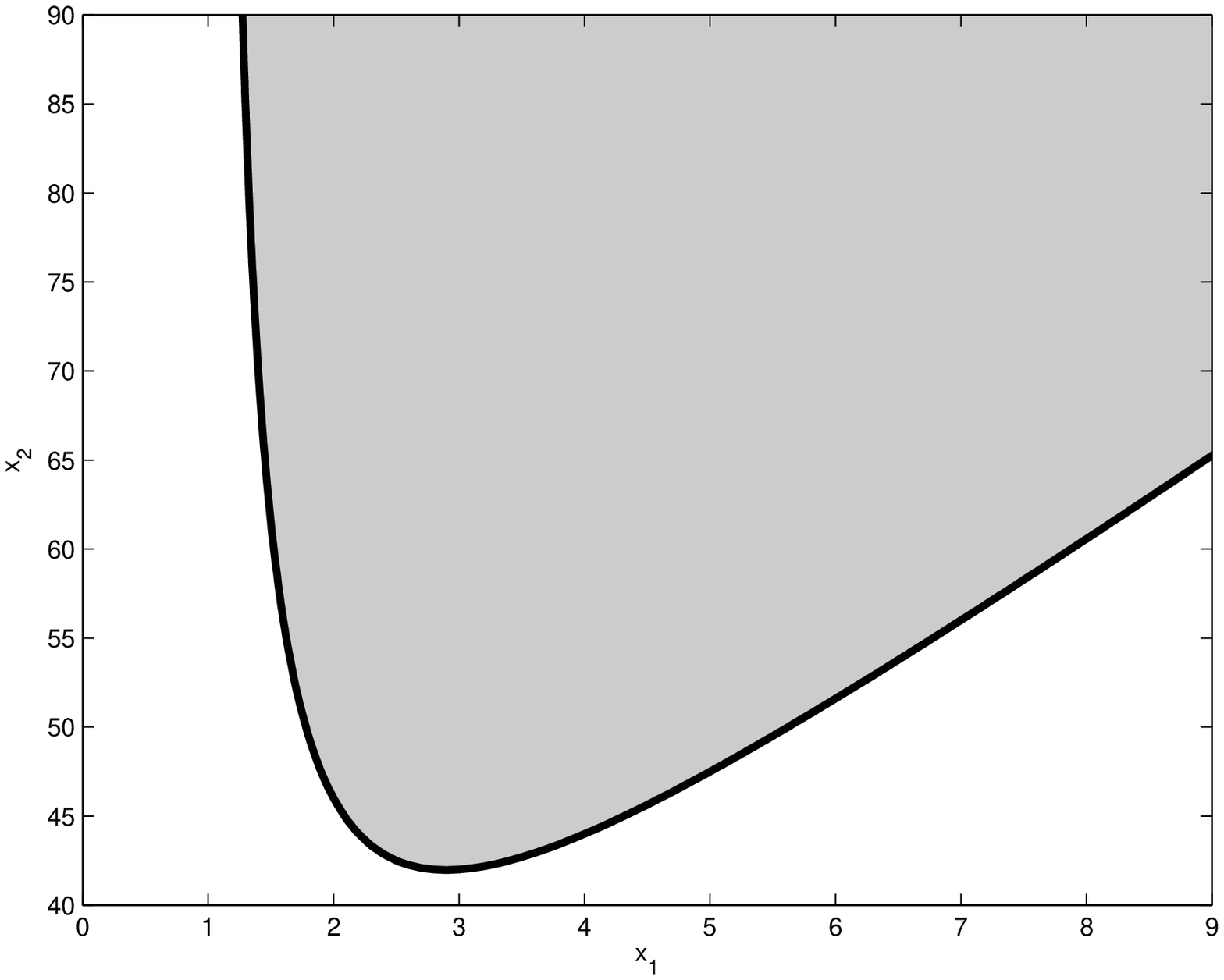}
\caption{{\tt NN1}.\label{fig-nn1}}
\end{center}
\end{minipage}
\hfill
\begin{minipage}[t]{.3\textwidth}
\begin{center}
\includegraphics[scale=0.3]{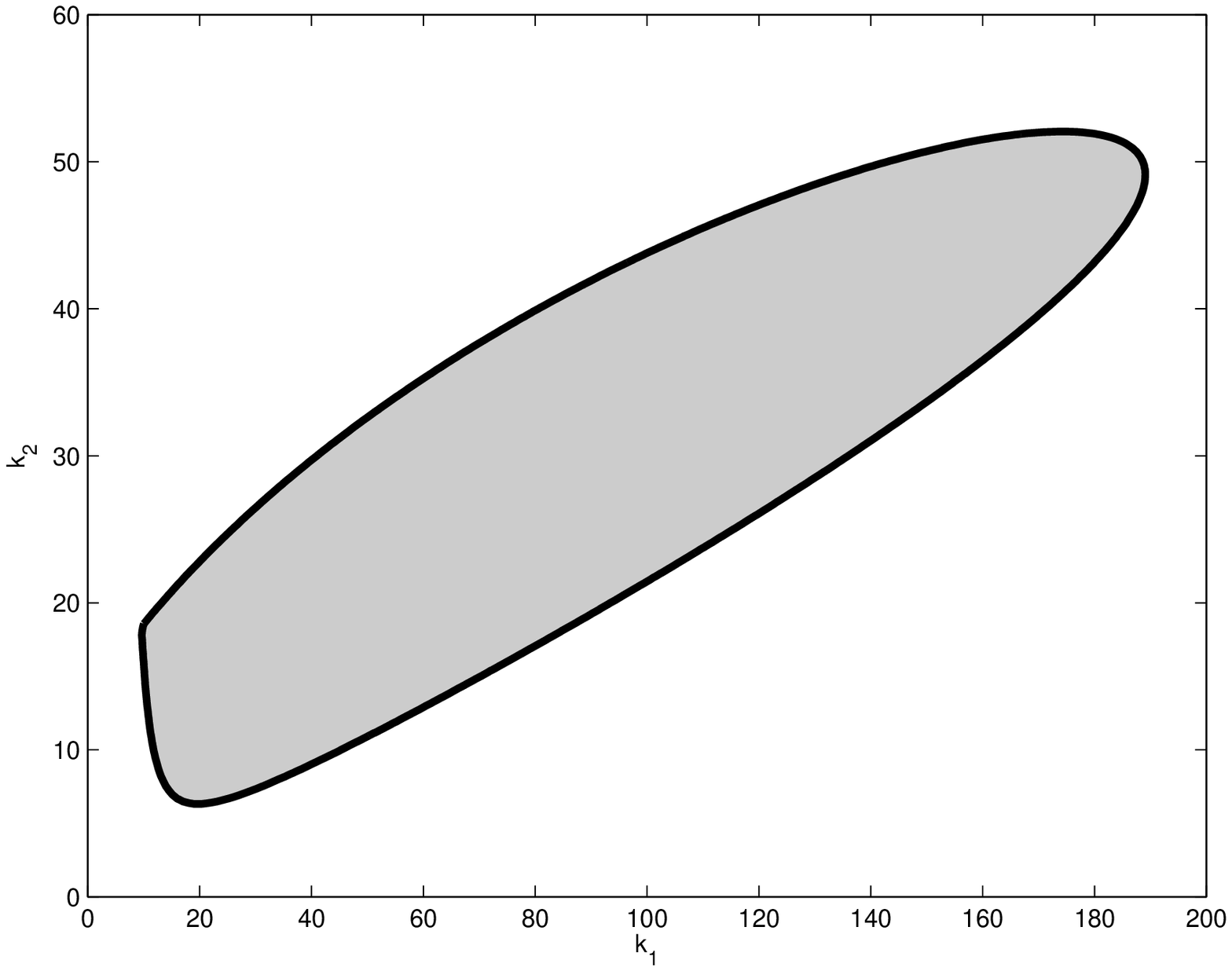}
\caption{{\tt NN5}.\label{fig-nn5}}
\end{center}
\end{minipage}
\hfill
\begin{minipage}[t]{.3\textwidth}
\begin{center}
\includegraphics[scale=0.3]{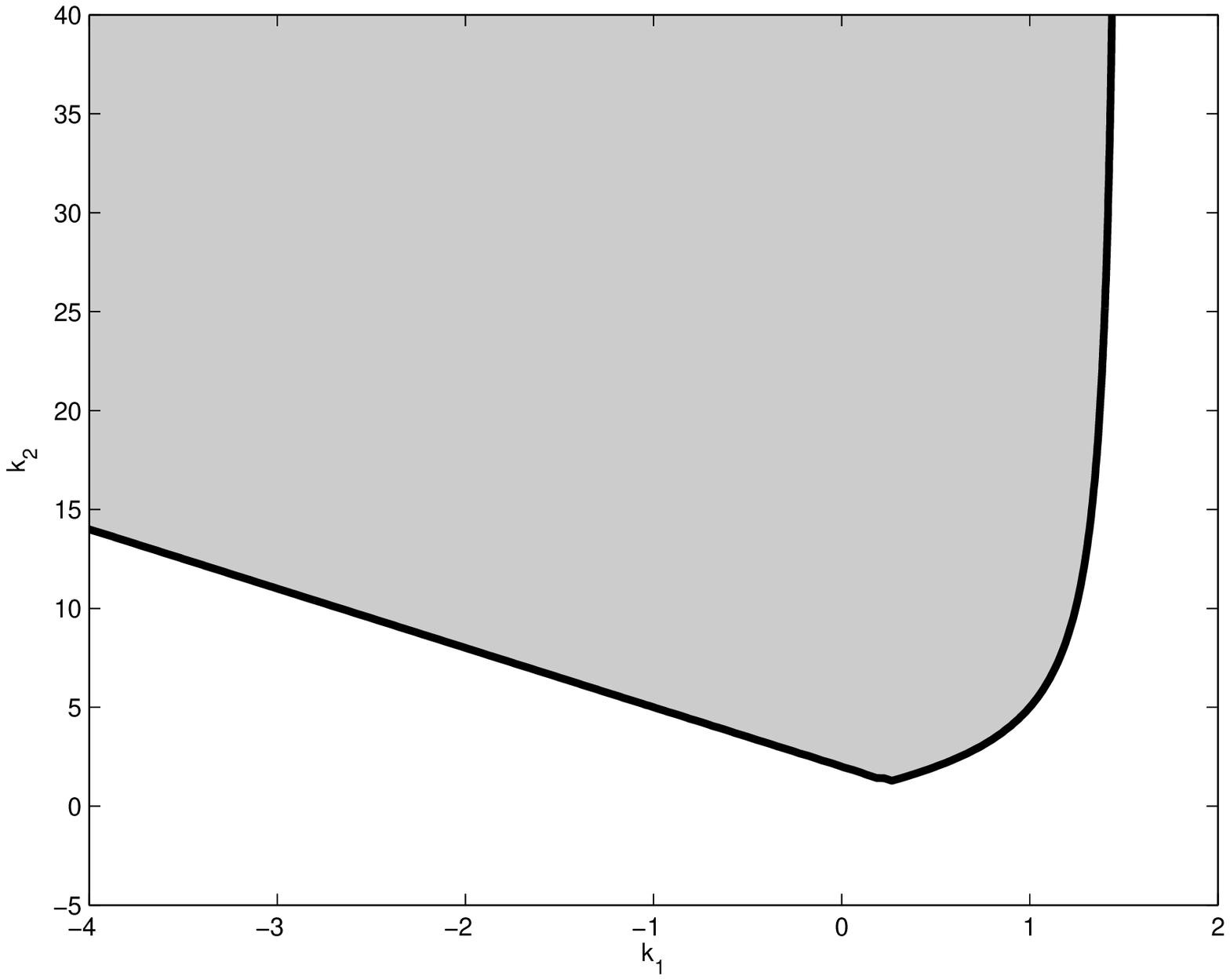}
\caption{{\tt NN17}.\label{fig-nn17}}
\end{center}
\end{minipage}
\end{figure}

\begin{figure}
\begin{center}
\includegraphics[scale=0.3]{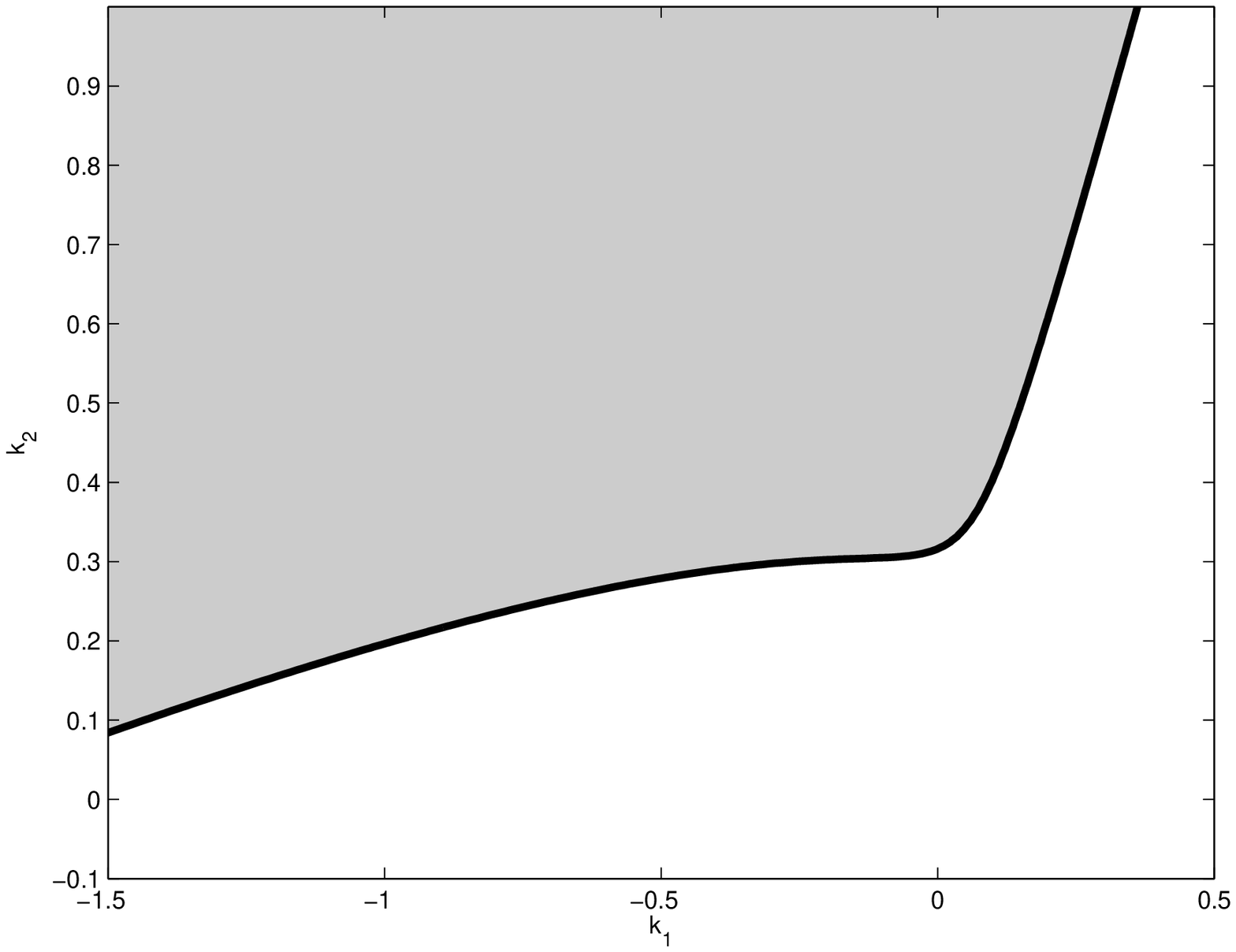}
\caption{{\tt HE1}.\label{fig-he1}}
\end{center}
\end{figure}

Consider the 7 two-parameter SOF
problems found in the database COMPleib \cite{compleib},
labelled {\tt AC4}, {\tt AC7}, {\tt AC17}, {\tt NN1}, {\tt NN5}, 
{\tt NN17} and {\tt HE1}. 
Stability regions are represented as shaded gray areas
on Figures \ref{fig-ac4} to
\ref{fig-he1}. Visual inspection reveals that
6 out of 7 stability regions seem to be convex. The only
apparently nonconvex example is {\tt HE1}.

In the remainder of the paper we will explain
why such planar stability regions are likely to be
convex, and how we can constructively derive
their LMI formulations when possible.

\section{Rational boundary of the stability region}

Define the curve
\[
{\mathcal C} = \{k \in {\mathbb R}^2 : p(j\omega,k) = 0, w \in {\mathbb R}\}
\]
which is the set of parameters $k$ for which polynomial
$p(s,k)$ has a root along the boundary of the stability
region, namely the imaginary axis.
Studying this curve is the key idea behind the
D-decomposition approach \cite{gryazina}. The curve
partitions the plane $(k_1,k_2)$ into regions
in which the number of stable roots of $p(s,k)$
remains constant. The union of regions for which
this number is equal to the degree of $p(s,k)$
is the stability region $\mathcal S$. Hence the
boundary of $\mathcal S$ is included in curve $\mathcal C$.

Note that $p(j\omega,k) = 0$ for some $w \in {\mathbb R}$
if and only if
\[
\begin{array}{rclcl}
p_R(\omega^2,k) & = & 
p_{0R}(\omega^2) + k_1 p_{1R}(\omega^2) + k_2 p_{2R}(\omega^2) & = & 0 \\
\omega p_I(\omega^2,k) & = &
\omega p_{0I}(\omega) + k_1 \omega p_{1I}(\omega) +
k_2 \omega p_{2I}(\omega) & = & 0\\
\end{array}
\]

Recall that we denote by $r_{\omega}(q_1,q_2)$ the
resultant of polynomials $q_1(\omega)$,$q_2(\omega)$
obtained by eliminating the scalar indeterminate $\omega$.
From the definition of the Hermite matrix, it holds
\begin{equation}\label{det}
h(k) = r_{\omega}(p_R(\omega^2,k),\omega p_I(\omega^2,k))
= \det H(k)
\end{equation}
from which the implicit algebraic description
\[
{\mathcal C} = \{k : h(k) = 0\}
\]
follows.

\begin{lemma}\label{factor}
The determinant of the Hermite matrix can be factored as
\[
h(k) = l(k)g(k)^2
\]
where $l(k)$ is affine, and $g(k)$ is
a generically irreducible polynomial.
\end{lemma}

{\bf Proof}:
The result follows from basic properties of resultants:
$h(k) = r_{\omega}(p_R(\omega^2,k),\omega p_I(\omega^2,k)) =
r_{\omega}(p_R(\omega^2,k),\omega^2) r_{\omega}(p_R(\omega^2,k),
p_I(\omega^2,k)) = r_{\omega}(p_R(\omega,k),\omega)
r_{\omega}(p_R(\omega,k),p_I(\omega,k))^2$. Take $g(k) =
r_{\omega}(p_R(\omega,k),p_I(\omega,k))$. Since $p_R(\omega,k)$
is affine in $k$ it follows that $l(k) = r_{\omega}(p_R(\omega,k),
\omega)$ is affine in $k$.$\Box$

The curve can therefore be decomposed as the union of a line
and a simpler algebraic curve
\[
{\mathcal C} = {\mathcal L} \cup {\mathcal G} = \{k : l(k) = 0\} \cup \{k : g(k) = 0\}.
\]
The equation of line $\mathcal L$ was already given in the proof
of Lemma \ref{factor}, namely
\[
l(k) = r_{\omega}(p_R(\omega), \omega) =
p_R(0,k) = p_{0R}(0) + k_1 p_{1R}(0) + k_2 p_{2R}(0).
\]
The defining polynomial of the other curve component $\mathcal G$ can be
obtained via the formula 
\[
g(k) = r_{\omega}(p_R(\omega,k),p_I(\omega,k)).
\]
From the relations
\[
\left[\begin{array}{cc}
p_{1R}(\omega^2) & p_{2R}(\omega^2) \\
\omega p_{1I}(\omega^2) & \omega p_{2I}(\omega^2) 
\end{array}\right]
\left[\begin{array}{c}
k_1 \\ k_2
\end{array}\right]  =
-\left[\begin{array}{c}
p_{0R}(\omega^2) \\ \omega p_{0I}(\omega^2)
\end{array}\right]
\]
we derive a rational parametrization of $\mathcal G$:
\begin{equation}\label{param}
\left[\begin{array}{c}
k_1(\omega^2) \\ k_2(\omega^2)
\end{array}\right] = \frac{
\left[\begin{array}{cc}
p_{2I}(\omega^2) & -p_{2R}(\omega^2) \\
-p_{1I}(\omega^2) & p_{1R}(\omega^2)
\end{array}\right]}
{p_{1I}(\omega^2)p_{2R}(\omega^2)
-p_{1R}(\omega^2)p_{2I}(\omega^2)}
\left[\begin{array}{c}
p_{0R}(\omega^2) \\ p_{0I}(\omega^2)
\end{array}\right]
= 
\left[\begin{array}{c}
\frac{q_1(\omega^2)}{q_0(\omega^2)} \\
\frac{q_2(\omega^2)}{q_0(\omega^2)} \\
\end{array}\right]
\end{equation}
which is well-defined since by assumption $p_1(s)/p_2(s)$
is not a constant. From this parametrization we can derive a
symmetric linear determinantal form of
the implicit equation of this curve.

\begin{lemma}\label{pencil}
The symmetric affine pencil
\[
G(k) = B_{\omega}(q_1,q_2)+k_1B_{\omega}(q_2,q_0)+k_2B_{\omega}(q_1,q_0).
\]
is such that
${\mathcal  G} = \{k : \det G(k) = 0\}$.
\end{lemma}

{\bf Proof:} Rewrite the system of equations (\ref{param}) as
\[
\begin{array}{rcccl}
a(\omega^2,k) & = & q_1(\omega^2) - k_1 q_0(\omega^2) & = & 0 \\
b(\omega^2,k) & = & q_2(\omega^2) - k_2 q_0(\omega^2) & = & 0 \\
\end{array}
\]
and use the B\'ezoutian resultant
to eliminate indeterminate $\omega$ and obtain conditions
for a point $(k_1,k_2)$ to belong to the curve. The B\'ezoutian
matrix is $B_{\omega}(a,b) = B_{\omega}(q_1-k_1q_0,q_2-k_2q_0) =
B_{\omega}(q_1,q_2)+k_1B_{\omega}(q_2,q_0)+k_2B(q_1,q_0)$.
Linearity in $k$ follows from bilinearity of the B\'ezoutian
and the common factor $q_0$.$\Box$

Finally, let $C(k) = \mathrm{diag}\:\{l(k), G(k)\}$
so that curve $\mathcal C$ can be
described as a determinantal locus
\[
{\mathcal C} = \{k \: :\: \det C(k) = 0\}.
\]

\section{LMI formulation}

Curve $\mathcal C$ partitions the plane into
several connected components, that we denote by ${\mathcal S}_i$
for $i=1,\ldots,N$.

\begin{lemma}\label{sc}
If $C(k) \succ 0$ for some point $k$ in the interior
of ${\mathcal S}_i$ for some $i$
then ${\mathcal S}_i = \{k : C(k) \succeq 0\}$
is a convex LMI region.
\end{lemma}

{\bf Proof:} 
Follows readily from the affine dependence of $C(k)$ on $k$
and from the fact that the boundary of ${\mathcal S}_i$
is included in $\mathcal C$.
$\Box$

Convex sets which
admit an LMI representation are called rigidly convex
in \cite{hv}. Rigid convexity is stronger
than convexity. It may happen that ${\mathcal S}_i$ is convex
for some $i$, yet $C(k)$ is not positive definite for
points $k$ within ${\mathcal S}_i$. 

\begin{lemma}\label{sh}
Stability region $\mathcal S$ is the union of sets
${\mathcal S}_i$ containing points $k$ such that
$H(k) \succ 0$.
\end{lemma}

{\bf Proof:}
Follows readily from Lemma \ref{hermite}.$\Box$

Note that it may happen that ${\mathcal S}_i$ is
convex LMI for some $i$, yet $H(k)$ is not
positive definite for points $k$ within ${\mathcal S}_i$.

\begin{corollary}
If $H(k) \succ 0$ and $C(k) \succ 0$ for some
point $k$ in the interior of ${\mathcal S}_i$,
then ${\mathcal S}_i$ is an LMI region included
in the stability region.
\end{corollary}

{\bf Proof:}
Combine Lemmas \ref{sc} and \ref{sh}.$\Box$

Quite often, on practical instances, we observe
that ${\mathcal S} = {\mathcal S}_i$ for some $i$
is a convex LMI region.

Practically speaking, once curve $\mathcal C$
is expressed as a determinantal locus, the search
of points $k$ such that $C(k) \succ 0$ can
be formulated as an eigenvalue problem,
but this is out of the scope of this paper.

\section{Examples}

\subsection{Example 1}

\begin{figure}[h]
\begin{center}
\includegraphics[scale=0.5]{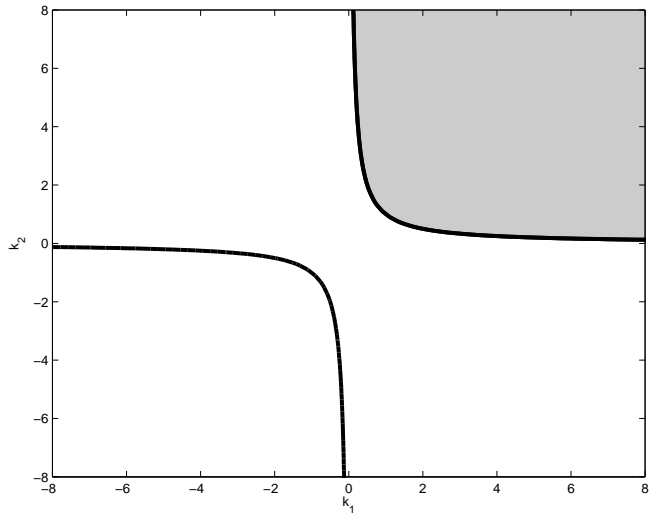}
\caption{Vishnegradsky's degree 3 polynomial.
Rational curve $\mathcal C$ 
with convex LMI stability region $\mathcal S$
(in gray)\label{fig-degree3}}
\end{center}
\end{figure}

As mentioned in \cite{gryazina}, Vishnegradsky in 1876
considered the polynomial $p(s,k) = s^3+k_1s^2+k_2s+1$
and concluded that its stability region
${\mathcal S} = \{k : k_1>0, k_1k_2>1\}$
is convex hyperbolic.

The Hermite matrix of $p(s,k)$ is given by
\[
H(k) = \left[\begin{array}{ccc} k_2 & 0 & 1 \\
0 & k_1 k_2 - 1 & 0\\
1 & 0 & k_1
\end{array}\right]
\]
and hence after a row and column permutation
the quadratic matrix inequality formulation
\[
{\mathcal S} = \{k : \left[\begin{array}{ccc}
k_2 & 1 & 0\\
1 & k_1 & 0 \\
0 & 0 & k_1 k_2 - 1
\end{array}\right] \succ 0\}
\]
explains why the region is convex. Indeed,
the determinant of the 2-by-2 upper matrix,
affine in $k$, is equal to the remaining
diagonal entry, which is here redundant.
The stability region can therefore be modeled as
the LMI
\[
{\mathcal S} = \{k : \left[\begin{array}{cc}
k_1 & 1 \\ 1 & k_2
\end{array}\right] \succ 0\}
\]
see Figure \ref{fig-degree3}.

\subsection{Example 2}

\begin{figure}[h]
\begin{center}
\includegraphics[scale=0.5]{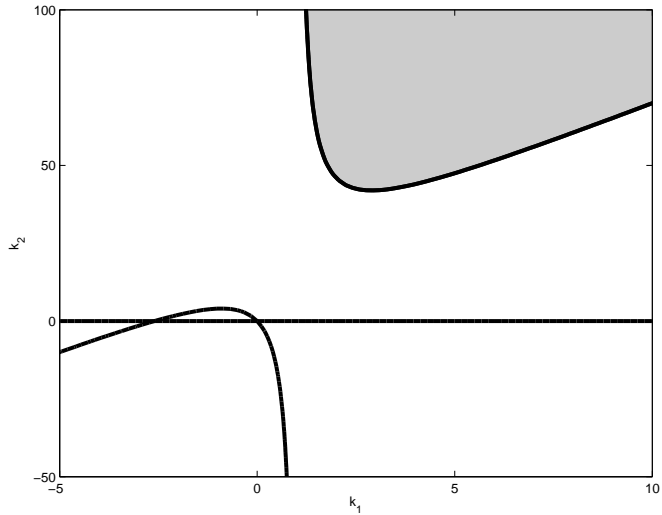}
\caption{{\tt NN1} SOF example.
Rational curve $\mathcal C$ 
with convex LMI stability region $\mathcal S$
(in gray)\label{fig-nn1b}}
\end{center}
\end{figure}

Consider problem {\tt NN1}
from \cite{compleib}, for which $p_0(s) = s(s^2-13)$,
$p_1(s) = s(s-5)$, $p_2(s) = s+1$ in (\ref{poly}).
We have $p_R(\omega^2) = -k_1\omega^2+k_2$ and $p_I(\omega^2) = -\omega^2
-13-5k_1+k_2$ and
\[
H(k) =
\left[\begin{array}{ccc}
k_2(-13-5k_1+k_2) & 0 & -k_2 \\
0 & k_1(-13-5k_1+k_2)-k_2 & 0 \\
-k_2 & 0 & k_1
\end{array}\right].
\]
Hence $h(k) = \det H(k) = k_2(-13k_1-k_2-5k_1^2+k_1k_2)^2$
and then $l(k)=k_2$, $g(k)=-13k_1-k_2-5k_1^2+k_1k_2$.
A rational parametrization of the curve ${\mathcal G} = \{k :
g(k)=0\}$ is given by
\[
\begin{array}{rcl}
k_1(\omega^2) & = & (\omega^2+13)/(\omega^2-5) \\
k_2(\omega^2) & = & \omega^2(\omega^2+13)/(\omega^2-5)
\end{array}
\]
from which we derive the symmetric affine determinantal 
representation ${\mathcal G} = \{k :
\det G(k)=0\}$ with
\[
G(k) = \left[\begin{array}{cc}
169+65k_1-18k_2 & 13+5k_1\\
13+5k_1 & 1-k_1
\end{array}\right].
\]
The pencil representing $\mathcal C$ is therefore
\[
C(k) = \left[\begin{array}{ccc}
k_2 & 0 & 0\\
0 & 169+65k_1-18k_2 & 13+5k_1 \\
0 & 13+5k_1 & 1-k_1
\end{array}\right].
\]

We can check that ${\mathcal S} = \{k : C(k) \succ 0\}$
is a convex LMI formulation of the stability region
represented on Figure \ref{fig-nn1}. Compare with
Figure \ref{fig-nn1b} where we represent also the
curve ${\mathcal C} = \{k : \det C(k) = 0\}$.

\subsection{Example 3}

\begin{figure}[h]
\begin{center}
\includegraphics[scale=0.5]{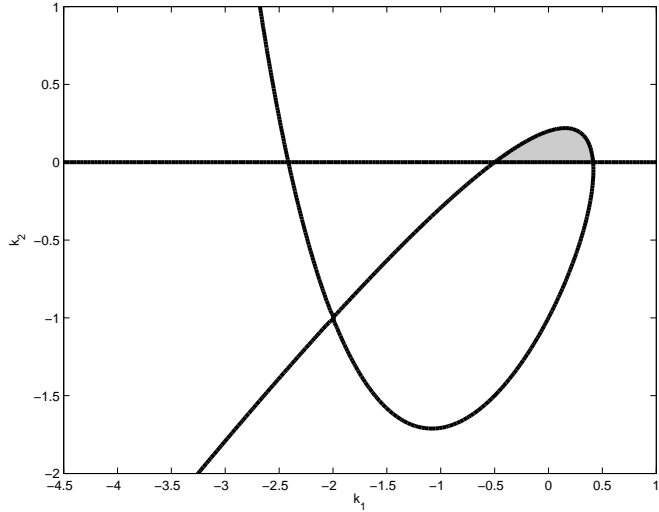}
\caption{Francis' example.
Rational curve $\mathcal C$ 
with convex LMI stability region $\mathcal S$
(in gray)\label{fig-francis87}}
\end{center}
\end{figure}

This example, originally from Francis (1987), is also described
in \cite{gryazina}. A SISO plant $(s-1)(s-2)/(s+1)(s^2+s+1)$
must be stabilized with a PI controller $k_1+k_2/s$.
Equivalently, $p_0(s) = s(s+1)(s^2+s+1)$, $p_1(s) = s(s-1)(s-2)$
and $p_2(s) = (s-1)(s-2)$ in (\ref{poly}).

For these values we obtain
\[
C(k) = 
\left[\begin{array}{cccc}
2k_2 & 0 & 0 & 0  \\
0 & 14+28k_1-54k_2 & -20-40k_1+18k_2 & 2+4k_1 \\
0 & -20-40k_1+18k_2 & 77-53k_1+36k_2 & -11+5k_1 \\
0 & 2+4k_1 & -11+5k_1 & 5+k_1
\end{array}\right]
\]
and the LMI stability region ${\mathcal S} = \{k : C(k) \succ 0\}$
represented on Figure \ref{fig-francis87} together with
the quartic curve ${\mathcal C} = \{k : \det C(k) = 0\}$.

\subsection{Example 4}

Consider \cite[Example 14.4]{ackermann}
for which $p_0(s) = s^4+2s^3+10s^2+10s+14+2a$,
$p_1(s) = 2s^3+2s-3/10$, $p_2(s) = 2s+1$, with
$a \in {\mathbb R}$ a parameter.

\begin{figure}[h]
\begin{center}
\includegraphics[scale=0.5]{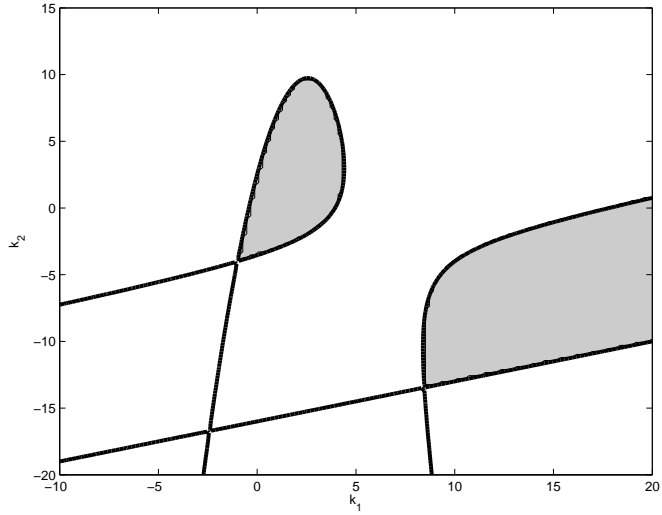}
\caption{Ackermann's example
with $a=1$. Rational curve $\mathcal C$ and
stability region $\mathcal S$ consisting of two
disconnected regions (in gray), the one
including the origin being LMI.
\label{fig-ackermann1}}
\end{center}
\end{figure}

We obtain $C(k) = \mathrm{diag}\:\{l(k), G(k)\}$ with
\[
\begin{array}{rcl}
l(k) & = & 140+20a-3k_1+10k_2 \\[1em]
G(k) & = & 
\left[\begin{array}{ccc}
\begin{array}{c}
7920+4860a+400a^2\\
+(-1609-60a)k_1+(-270+200a)k_2
\end{array} & \star & \star \\
-8350-2000a+1430k_1+130k_2 & 8370-1230k_1-100k_2 & \star\\
900+200a-130k_1 & -900+100k_1 & 100
\end{array}\right]
\end{array}
\]
where symmetric entries are denoted by stars.

When $a=1$, the stability region consists of
two disconnected components. The one
containing the origin $k_1=k_2=0$ is the
LMI region $\{k : C(k) \succ 0\}$,
see Figure \ref{fig-ackermann1}.

\begin{figure}[h]
\begin{center}
\includegraphics[scale=0.5]{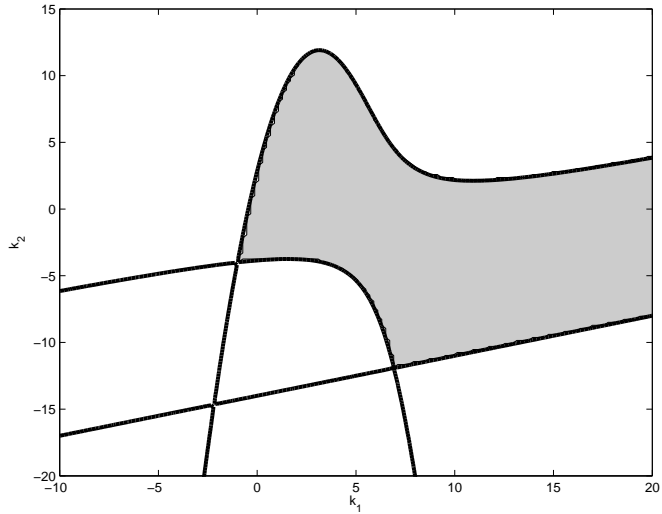}
\caption{Ackermann's example
with $a=0$.
Rational curve $\mathcal C$ and non-convex
stability region $\mathcal S$ (in gray).
\label{fig-ackermann2}}
\end{center}
\end{figure}

When $a=0$, the stability region $\mathcal S$
is the non-convex region represented on Figure
\ref{fig-ackermann2}. The LMI region 
$\{k : C(k) \succ 0\}$ is not included
in $\mathcal S$ in this case.

\section{Conclusion}

In this paper we have explained why the planar stability
region of a polynomial may be convex with
an explicit LMI representation. This is an instance
of hidden convexity of a set which is otherwise described
by intersecting (generally non-convex) Routh-Hurwitz
minors sublevel sets or by enforcing positive definiteness
of a (generally non-convex) quadratic Hermite matrix.

Practically speaking,
optimizing a closed-loop performance criterion
over an LMI formulation of the stability region
is much simpler than optimizing over the non-linear
formulation stemming from the Routh-Hurwitz minors
or the Hermite quadratic matrix inequality.

Convexity in the parameter space was already exploited
in \cite{ho,ak} in the context of PID controller design.
It was shown that when the proportional gain is
fixed, the set of integral and derivative gains
is a union of a finite number of polytopes.

Extension of these ideas to the case of more than 2 parameters
seems to be difficult. The problem of finding a
symmetric affine determinantal representation
of rationally parametrized surfaces or hypersurfaces
is not yet well understood, to the best of our knowledge.
For example, in the simplest third
degree case $p(s,k) = s^3+k_1s^2+k_2s+k_3$, how could
we find four symmetric real matrices $A_0$, $A_1$, $A_2$, $A_3$ satisfying
$\det (A_0 + A_1 k_1 + A_2 k_2 + A_3 k_3) = k_1k_3-k_2$ ?

\end{document}